\documentclass[3p]{elsarticle}
\usepackage[all]{xy}
\usepackage{color}
\usepackage{url}
\usepackage{amscd, amssymb, amsfonts, amsthm, amsmath}
\usepackage{mathdots}
\usepackage{hyperref}
\usepackage{multicol}

\newcommand{\topp}{\mbox{\rm{Top}}}
\newcommand{\length}{\mbox{\rm{l}}}
\newcommand{\gl}{\mbox{\rm{gldim}}}

\newcommand{\fidim}{\mbox{\rm{$\phi$dim}}}

\newcommand{\findim}{\mbox{\rm{findim}}}

\def\mod{\mbox{\rm{mod}}}
\def\Mod{\mbox{\rm{Mod}}}
\def\ind{\mbox{\rm{ind}}}
\def\add{\mbox{\rm{add}}}

\def\pd{\mbox{\rm{pd}}}
\def\id{\mbox{\rm{id}}}

\def\enn{\hbox{\rm{End}}}
\def\rk{\hbox{\rm{rk}}}
\def\rad{\hbox{\rm{rad}}}
\def\ll{\hbox{\rm{ll}}}
\def\supp{\hbox{\rm{supp}}}

\def\Orb{\mbox{\rm{Orb}}}

\def\Rep{\mbox{\rm{Rep}}}
\def\rep{\mbox{\rm{rep}}}

\def\start{\mbox{\rm{s}}}
\def\target{\mbox{\rm{t}}}

\begin{document}
\newcommand{\mono}[1]{%
\gdef\puA{#1}}
\newcommand{\puA}{}
\newcommand{\faculty}[1]{%
\gdef\puC{#1}}
\newcommand{\puC}{}
\newcommand{\facultad}[1]{%
\gdef\puD{#1}}
\newcommand{\puD}{}
\newcommand{\N}{\mathbb{N}}
\newcommand{\Z}{\mathbb{Z}}
\newtheorem{teo}{Theorem}[section]
\newtheorem{prop}[teo]{Proposition}
\newtheorem{lema}[teo] {Lemma}
\newdefinition{ej}[teo]{Example}
\newtheorem{obs}[teo]{Remark}
\newtheorem{defi}[teo]{Definition}
\newtheorem{coro}[teo]{Corollary}
\newtheorem{nota}[teo]{Notation}



\title{The Igusa-Todorov $\phi$-dimension on Morita context algebras}

\author[add]{Marcos Barrios}
\ead{marcosb@fing.edu.uy} 

\author[add]{Gustavo Mata\corref{cor}}
\ead{gmata@fing.edu.uy}

\address[add]{Universidad de La Rep\'ublica, Facultad de Ingenier\'ia -  Av. Julio Herrera y Reissig 565, Montevideo, Uruguay}
\cortext[cor]{Corresponding Author}

\begin{abstract}
In this article we prove that, under certain hypotheses, Morita context algebras that have zero bimodule morphisms have finite $\phi$-dimension. We also study the behaviour of the $\phi$-dimension for an algebra and its opposite. In particular we show that the $\phi$-dimension of an Artin algebra is not symmetric, i.e. there exists a finite dimensional algebra $A$ such that $\fidim (A) \not = \fidim (A^{op})$. 
\end{abstract} 


\maketitle

\section{Introduction}

For an Artin algebra $A$ the finitistic dimension is defined as 
$$\findim(A) = \sup\{ \pd(M) |\ M \in \mod A,\ \pd M < \infty \}.$$
The (small) finitistic dimension conjecture states that $\findim(A) < \infty$.
In an attept to prove this conjecture, K. Igusa and G. Todorov introduced the functions $\phi$ and $\psi$ in \cite{IT}, nowadays called the Igusa-Todorov functions.
Later, using the Igusa-Todorov functions, F. Huard and M. Lanzilotta introduced the $\phi$-dimension and the $\psi$-dimension in \cite{HL}, and proved that an Artin algebra is selfinjective if and only if its $\phi$-dimension (or its $\psi$-dimension) is zero. Hence these dimensions show how far an Artin algebras is from being selfinjective, in some sense.

Since M. Auslander proved that the left global dimension is equal to the right global dimension for semi-primary rings,  the following question arises naturally due to the relationship of the $\phi$-dimension with the global dimension: Is the $\phi$-dimension symmetric? i.e. for an Artin algebra $A$, $\phi \dim(A) = \phi \dim(A^{op})$?
In the case of the $\psi$-dimension it is not true and an example of radical square zero algebras is given in \cite{LMM}.   

The symmetry conjecture of the $\phi$-dimension was proved in \cite{LMM} for radical square zero algebras, in \cite{LM} for Gorenstein algebras, and in \cite{BMR} for truncated path algebras. On the other hand the conjecture is not true for semiperfect coalgebras. In \cite{HLM17} was given an easy example (Example 2.4) of a semiperfect coalgebra $C$ such that $\fidim (C) \neq \fidim(C^{op})$.

Let $A$ and $B$ be two Artin algebras, $Y$ an $A$-$B$-bimodule, $X$ a $B$-$A$-bimodule,
$\alpha : X \otimes_A Y \rightarrow B$ a $B$-$B$-bimodule homomorphism, and $\beta : Y \otimes_B X \rightarrow A$ an $A$-$A$-bimodule homomorphism. Then from the Morita context $M = (A, Y, X, B, \alpha, \beta)$ we define the Morita context algebra:
$$ \Lambda_{(\alpha, \beta)} = \begin{pmatrix}
A & Y\\
X & B
\end{pmatrix},$$
where the multiplication in $\Lambda_{(\alpha, \beta)}$ is given by 
$$\begin{pmatrix}
a & y\\
x & b
\end{pmatrix} \cdot \begin{pmatrix}
a' & y'\\
x' & b'
\end{pmatrix} = \begin{pmatrix}
aa'+\beta(y\otimes x') & ay'+yb'\\
xa' + bx' & bb' + \alpha(x \otimes y')
\end{pmatrix}$$
and the maps $\alpha$ and $\beta$ satisfy $\alpha(x \otimes y)x' = x\beta(y \otimes x' )$ and $y\alpha(x \otimes y') = \beta(y \otimes x)y' $, to make $\Lambda_{(\alpha, \beta)}$ associative.

We say that a Morita context algebra has zero bimodule homomorphisms if $\alpha = \beta = 0$.

The family of Morita context algebras with zero bimodule homomorphisms has been extensively studied from the homological point of view. We can cite for example, \cite{CRS}, \cite{GP*} and \cite{GP}.  \\

In \cite{BM1} the authors prove that, under some hypotheses on the bimodule $Y$, a triangular algebra 
$$\begin{pmatrix}
A & Y\\
0 & B
\end{pmatrix}$$ has finite $\phi$-dimension if $A$ and $B$ have finite $\phi$-dimensions. 
A natural extension of triangular matrix algebras is the class of Morita context algebras. For that reason, in this article we focus on Morita context algebras that have the folowing shape: 
If $A = \frac{\Bbbk Q_A}{I_A}$ and $B = \frac{\Bbbk Q_B}{I_B}$ are finite dimensional algebras, $C = \frac{\Bbbk \Gamma}{I_C}$ has the following hypotheses

\begin{itemize}

\item[{\bf H1}:] ${Q_C}_0 = {Q_A}_0 \cup {Q_B}_0$.

\item[{\bf H2}:] ${Q_C}_1 = {Q_A}_1 \cup {Q_B}_1 \cup \{\alpha_j : \start(\alpha_j) \in {Q_A}_0, \ \target(\alpha_j) \in {Q_B}_0 \}_{j \in J} \cup  \{\beta_k : \start(\beta_k) \in {Q_B}_0, \ \target(\beta_k) \in {Q_A}_0 \}_{k \in K}$.

\item[{\bf H3}:] $\langle I_A, I_B, \alpha \alpha_j, \beta \beta_k \text{ for }\alpha \in {Q_A}_1, \beta\in {Q_B}_1, \alpha_j\beta_k, \beta_k\alpha_j \text{ where } j\in J, k\in K\rangle \subset I_C$.

\item[{\bf H4}:] $\mathcal{O}= [\add \Orb_{\Omega_A}(\Pi_A(\Omega_C(B_0))) \times \add \Orb_{\Omega_B}( \Pi_B(\Omega_C(A_0)))] \subset K_0(C)$ is syzygy finite.

\end{itemize}

Note that: 

\begin{itemize}

\item Hypotheses {\bf H1}, {\bf H2} and {\bf H3} imply that for all $M \in \mod C$, $\Omega(M) = N_A \oplus N_B$ where $Y_A \in \mod A$ and $Y_B \in \mod B$. They also imply that the bimodules are left-semisiple.

\item Hypotheses {\bf H1}, {\bf H2} and {\bf H3} with an equality in Hypothesis {\bf H3} imply that the bimodules $X$ and $Y$ are right-projective.

\item {\bf H4} is straightforward if $A$ and $B$ are syzygy finite or in case of the inclusion in 
Hypothesis {\bf H3} is an equality.

\end{itemize}

The content of this article can be summarised as follows.
Section 2 is devoted to collect some necessary material for the developing of this work.
In Section 3, we prove that Morita context algebras with hypotheses {\bf H1}, {\bf H2}, {\bf H3} and {\bf H4} which come from algebras with finite $\phi$-dimension also has finite $\phi$-dimension (Theorem 3.2). The previous result is a generalization of Theorem 5.2 of \cite{BM1}. We also give an explicit bound when the inclusion of the ideal $I_C$, in Hypothesis 3, is an equality. Similar results are obtained for the finitistic dimension. 
Finally, in Section 4, we exhibit an example which shows a finite dimensional algebra $A$ such that $\fidim (A) \not = \fidim (A^{op})$.

\section{Preliminaries}

Throughout this article $A$ is an Artin algebra and $\mod A$ is the category of finitely generated right $A$-modules, $\ind A$ is the subcategory of $\mod A$ formed by all indecomposable modules, $\mathcal{P}_A \subset \mod A$ is the class of projective $A$-modules. $\mathcal{S} (A)$ is the set of isoclasses of simple $A$-modules and $A_0 = \oplus_{S \in \mathcal{S}(A)}S$. For $M\in \mod A$ we denote by $M^k = \oplus_{i=1}^k M$, by $P(M)$ its projective cover, by $\Omega(M)$ its syzygy and by $\ll(M)$ its Loewy length. The set $\Orb_{\Omega(M)}$ is the $\Omega$-orbit of the module $M$, i.e. $\Orb_{\Omega(M)} = \{\Omega^n(M)\}_{n\geq 0}$. For a subcategory $\mathcal{C} \subset \mod A$, we denote by $\findim (\mathcal{C})$, $\gl (\mathcal{C})$ its finitistic dimension and its global dimension respectively and by $\add \mathcal{C}$ the full subcategory of $\mod A$ formed by all the sums of direct summands of every $M \in \mathcal{C}$. 

If $Q = (Q_0,Q_1,\start,\target)$ is a finite connected quiver, $\mathfrak{M}_Q$ denotes its adjacency matrix and $\Bbbk Q$ its associated {\bf path algebra}. We compose paths in $Q$ from left to right. 
Given $\rho$ a path in $\Bbbk Q$, $\length(\rho)$, $\start(\rho)$ and $\target(\rho )$ denote the length, start and target of $\rho$ respectively. For a quiver $Q$, we denote by $J_Q$ the ideal of $\Bbbk Q$ generated by all its arrows. If $I$ is an admissible ideal of $\Bbbk Q$ ($J^2_Q \subset I \subset J^m_Q$ for some $m \geq 2$), we say that $(Q,I)$ is a {\bf bounded quiver} and the quotient algebra $\frac{\Bbbk Q}{I}$ is the {\bf bound quiver algebra}. A {\bf relation} $\rho$ is an element in $\Bbbk Q$ such that $\rho = \sum \lambda_i w_i$ where the $\lambda_i$ are scalars (not all zero) and the $w_i$ are paths with $\length(w_i) \geq 2$ such that $\start(w_i) = \start(w_j)$ and $\target(w_i) = \target (w_j)$ if $i \not = j$. 
We recall that an admissible ideal $I$ is always generated by a finite set of relations (for a proof see Chapter II.2 Corollary 2.9 of \cite{ASS}).

For a quiver $Q$, we say that $M = (M^v, T_{\alpha})_{v\in Q_0, \alpha \in Q_1}$ is a {\bf representation} if
\begin{itemize}
\item $M^v$ is a $\Bbbk$-vectorial space for every $v\in Q_0$,
\item $T_{\alpha} :M^{\start(\alpha)} \rightarrow M^{\target(\alpha)}$ for every $\alpha \in Q_1$.
\end{itemize}

A representation $(M^v, T_{\alpha})_{v\in Q_0, \alpha \in Q_1}$ is finite dimensional if $M^v$ is finite dimensional for every $v \in Q_0$.

For a path $w = \alpha_1\ldots \alpha_n$ we define $T_{w} = T_{\alpha_1} \ldots T_{\alpha_n}$, and for a relation $\rho = \sum \lambda_i w_i $ we define $T_{\rho} = \sum \lambda_i T_{w_i}$. 

A representation $M = (M^v, T_{\alpha})_{v\in Q_0, \alpha \in Q_1}$ of $Q$ is {\bf bound by} $I$ if we have  $T_{\rho} = 0$ for all relations $\rho \in I$.

Let $M= (M^v, T_{\alpha})_{v\in Q_0, \alpha \in Q_1}$ and $\bar{M} = (\bar{M}^v, \bar{T}_{\alpha})_{v\in Q_0, \alpha \in Q_1}$ be two representations of the bounded quiver $(Q, I)$, a {\bf morphism} $f : M \rightarrow \bar{M}$ is a family $f=(f_v)_{v\in Q_0}$ of $\Bbbk$-linear maps $f_v: M^v \rightarrow \bar{M}^v$ such that for all arrow $\alpha : v \rightarrow w$ we have the following commutative diagram:

$$\xymatrix{ M^v \ar[r]^{T_{\alpha}} \ar[d]_{f_v} & M^w \ar[d]^{f_w}  \\ \bar{M}^v \ar[r]^{\bar{T}_{\alpha}}& \bar{M}^w}$$ 

We denote by $\Rep_{\Bbbk} (Q,I)$ the category of representations of $(Q, I)$ and by $\rep_{\Bbbk} (Q,I)$ the subcategory of $\Rep_{\Bbbk} (Q,I)$ formed by finite dimensional representations of $(Q, I)$. 

We recall that there is a $\Bbbk$-linear equivalence of categories $$F: \Rep_{\Bbbk} (Q,I) \rightarrow \Mod \frac{\Bbbk Q}{I}$$ that restricts to an equivalence of categories $G: \rep_{\Bbbk} (Q,I) \rightarrow \mod \frac{\Bbbk Q}{I}$ (for a proof see Chapter III.1, theorem 1.6 of \cite{ASS}).

If $A = \frac{\Bbbk Q_A}{I_A}$ and $B = \frac{\Bbbk Q_B}{I_B}$ are finite dimensional algebras, and $C = \frac{\Bbbk \Gamma}{I_C}$ the square algebra has the hypotheses {\bf H1} and {\bf H2}, then the functors
$\prod_A: \mod C \rightarrow \mod A$ and $\prod_B: \mod C \rightarrow \mod B$ are restictions of the representations, i.e. For a representation $M=(M^v, T_{\alpha})_{v\in {Q_C}_0, \alpha \in {Q_C}_1} \in \Rep_{\Bbbk}(Q_C,I_C)$, we have 

\begin{itemize}
\item $\prod_A(M) = (M^v, T_{\alpha})_{v\in {Q_A}_0, \alpha \in {Q_A}_1} $ and
\item  $\prod_B(M) = (M^v, T_{\alpha})_{v\in {Q_B}_0, \alpha \in {Q_B}_1}$. 
\end{itemize}
For a morphism $f = (f_v)_{v\in {Q_C}_0}: M\rightarrow \bar{M}$, then 

\begin{itemize}
\item $\prod_A(f) = (f_v)_{v\in {Q_A}_0}: \prod_A(M) \rightarrow \prod_A(\bar{M})$ and
\item  $\prod_B(f) = (f_v)_{v\in {Q_B}_0}: \prod_B(M) \rightarrow \prod_B(\bar{M})$.
\end{itemize}

\subsection{Igusa-Todorov functions}

In this section we exhibit some general facts about the Igusa-Todorov functions for an Artin algebra A.
The aim is to introduce material which will be used in the following sections.

\begin{lema}(Fitting Lemma)
Let $R$ be a noetherian ring. Consider a left $R$-module $M$ and $f \in \enn_R (M)$. Then, for any finitely generated $R$-submodule $X$ of $M$, there is a non-negative integer
$$\eta_{f}(X)= \min\{ k \text{ a non-negative integer}: f \vert _{f^m (X)} : f^m (X) \rightarrow f^{m+1} (X), \text{ is injective } \forall m \geq k\}.$$
Furthermore, for any $R$-submodule $Y$ of $X$, we have that $\eta_f(Y) \leq \eta_f(X)$.
\end{lema}

\begin{defi}(\cite{IT})
Let  $K_0(A)$ be the abelian group generated by all symbols $[M]$, with $M \in \mod A$, modulo the relations
\begin{enumerate}
  \item $[M]-[M']-[M'']$ if  $M \cong M' \oplus M''$,
  \item $[P]$ for each projective module $P$.
\end{enumerate}
\end{defi}

Let $\bar{\Omega}: K_0 (A) \rightarrow K_0 (A)$ be the group endomorphism induced by $\Omega$, i.e. $\bar{\Omega}([M]) = [\Omega(M)]$. We denote by $K_i (A) = \bar{\Omega}(K_{i-1}(A))= \ldots = \bar{\Omega}^{i}(K_{0} (A))$ for $i \geq 1$. 
We say that a subgroup $G \subset K_0(A)$ is syzygy finite if there is $n \geq 0$ such that $\bar{\Omega}^n(G)$ is finitely generated. 
For $M\in \mod A$, $\langle \add M\rangle$ denotes the subgroup of $K_0 (A)$ generated by the classes of indecomposable summands of $M$.

 For a subcategory $\mathcal{C} \subset \mod A$, we denote by $\langle\mathcal{C}\rangle \subset K_0(A)$ the free abelian group generated by the classes of direct summands of modules of $\mathcal{C}$.

\begin{defi}\label{monomorfismo}(\cite{IT}) 
The \textbf{(right) Igusa-Todorov function} $\phi$ of $M\in \mod A$  is defined as 
\[\phi_{A}(M) = \eta_{\bar{\Omega}}(\langle \add M \rangle).\]
 In case there is no possible misinterpretation we will use the notation $\phi$ for the Igusa-Todorov $\phi$ function.
\end{defi}

\begin{prop}(\cite{HLM}, \cite{IT}) \label{it1}
If $M,N\in\mod A$, then we have the following.

\begin{enumerate}
  \item $\phi(M) = \pd (M)$ if $\pd (M) < \infty$.
  \item $\phi(M) = 0$ if $M \in \ind A$ and $\pd(M) = \infty$.
  \item $\phi(M) \leq \phi(M \oplus N)$.
  \item $\phi\left(M^{k}\right) = \phi(M)$ for $k \geq 1$.
  \item $\phi(M) \leq \phi(\Omega(M))+1$.
\end{enumerate}

\end{prop}

\begin{prop}\label{invariante}
Suppose $G \subset K_0(A)$ is a finitely generated subgroup with $\rk(G) = m$.
\begin{enumerate}
\item If $\bar{\Omega}(G) \subset G$, then $\eta_{\bar{\Omega}}{}_{\vert G} \leq m$.

\item If $G$ is syzygy finite, then $\eta_{\bar{\Omega}}{}_{\vert G} < \infty$

\end{enumerate}
\begin{proof}
The proof of item 1 is similar to the proof of Proposition 3.6 (item 3) from \cite{LMM}.
The proof of item 2 is similar to the proof of Theorem 3.2. from  \cite{LMM}.
\end{proof}
\end{prop}

The result below follows directly from the fact that the Igusa-Todorov function verifies 
$$\phi(M) = \min\{l:\Omega{\vert}_{\Omega^{l+s} \langle \add M \rangle} \text{ is a monomorphism } \forall s \in \mathbb{N}\}.$$

\begin{prop}
Given $M \in  \mod A$, 
$$\phi(A) = \max\{n\in \mathbb{N}: \bar{\Omega}^n(v)=0\text{ and }\bar{\Omega}^{n-1}(v)\not=0 \text{ for some } v\in \langle \add M \rangle\}.$$

\end{prop}

\section{Morita context algebras and the Igusa-Todorov $\phi$ function}
In this section we compute the $\phi$-dimension of the Morita context algebra $C$ under certain hypotheses on $A$, $B$,  the quiver $Q_C$ and the ideal $I_C$.
 
The proof of the following lemma is similar to Theorem 5.2 of \cite{BM1}, for this reason it is left to the reader.

\begin{lema}\label{lemita}
Let $A = \frac{\Bbbk Q_A}{I_A}$ and $B = \frac{\Bbbk Q_B}{I_B}$ be finite dimensional algebras. Consider $C = \frac{\Bbbk \Gamma}{I_C}$ with the following conditions:

\begin{itemize}

\item ${Q_C}_0 = {Q_A}_0 \cup {Q_B}_0$.

\item ${Q_C}_1 = {Q_A}_1 \cup {Q_B}_1 \cup \{\alpha_j : \start(\alpha_j) \in {Q_A}_0, \ \target(\alpha_j) \in {Q_B}_0 \}_{j \in J} \cup  \{\beta_j : \start(\beta_k) \in {Q_B}_0,\ \target(\beta_k) \in {Q_A}_0 \}_{k \in K}$.

\item $\langle I_A, I_B, \alpha \alpha_j, \beta \beta_k \text{ for }\alpha \in {Q_A}_1, \beta\in {Q_B}_1, \alpha_j\beta_k, \beta_k\alpha_j \text{ where } j\in J, k\in K\rangle \subset I_C$.

Then for every module $M \in \mod C$, $\Omega_C(M) = M_1 \oplus M_2$ with $M_1 \in \mod A$ and $M_2 \in \mod B$. Suppose $\Omega_C(\topp(M))= \bar{M}_1 \oplus \bar{M}_2$ with $\bar{M}_1 \in \mod A$ and $\bar{M}_2 \in \mod B$, in case $M \in \mod A$, then $\bar{M}_2 = M_2$, in case $M \in \mod B$, then $\bar{M}_1 = M_1$. 
\end{itemize}

\end{lema}

\begin{teo} \label{teo grande}
Let $A_i = \frac{\Bbbk Q_{A_i}}{I_i}$ for $i=1,\ldots, m$ be finite dimensional algebras. Consider $C = \frac{\Bbbk \Gamma}{I_C}$ with the following conditions:

\begin{itemize}

\item ${Q_C}_0 = \bigcup_{i=1}^m {Q_{A_i}}_0$.

\item ${Q_C}_1 = \bigcup_{i=1}^m{Q_{A_i}}_1  \cup  \{\delta^{i,j}_k : \start(\delta^{i,j}_k ) \in {Q_{A_{i}}}_0, \ \target(\delta^{i,j}_k ) \in {Q_{A_j}}_0 \}_{k \in K_{i,j}, i \neq j, i,j \in \{1, \ldots, m\}} $.

\item $\langle I_i, \  \alpha \delta^{i,j}_k, \delta^{i,j}_k \delta^{i',j'}_{k'} \text{ for }\alpha \in {Q_{A_i}}_1, \text{ where } i,i',j,j' \in \{1, \ldots, m\}, k \in K_{i,j}, k'\in K_{i',j'}  \rangle \subset I_C$.

\item $\mathcal{O}= [\prod_{i=1}^m \add \Orb_{\Omega_A}(\prod_{A_i}(\Omega_C( \oplus_{j \not=i} ({A_j}_0) )))] \subset K_0(C)$ is syzygy finite.


\item $\fidim(A_i) < \infty, \forall i = 1, \ldots m$. 

\end{itemize}

Then $\fidim(C) < \infty$.

\end{teo}

The proof of the previous result is analogous to the proof of the following result.

\begin{teo}\label{teo}
Let $A = \frac{\Bbbk Q_A}{I_A}$ and $B = \frac{\Bbbk Q_B}{I_B}$ be finite dimensional algebras. Consider $C = \frac{\Bbbk \Gamma}{I_C}$ with the following conditions:

\begin{itemize}

\item ${Q_C}_0 = {Q_A}_0 \cup {Q_B}_0$

\item ${Q_C}_1 = {Q_A}_1 \cup {Q_B}_1 \cup \{\alpha_j : \start(\alpha_j) \in {Q_A}_0, \ \target(\alpha_j) \in {Q_B}_0 \}_{j \in J} \cup  \{\beta_k : \start(\beta_k) \in {Q_B}_0, \ \target(\beta_k) \in {Q_A}_0 \}_{k \in K}$

\item $\langle I_A, I_B, \alpha \alpha_j, \beta \beta_k \text{ for }\alpha \in {Q_A}_1, \beta\in {Q_B}_1, \alpha_j\beta_k, \beta_k\alpha_j \text{ where } j\in J, k\in K\rangle \subset I_C$.

\item $\mathcal{O}= [\add \Orb_{\Omega_A}(\Pi_A(\Omega_C(B_0))) \times \add \Orb_{\Omega_B}( \Pi_B(\Omega_C(A_0)))] \subset K_0(C)$ is syzygy finite.


\item $\fidim(A) < \infty$ and $\fidim(B) < \infty$.

\end{itemize}

Then $\fidim(C) < \infty$

\begin{proof}
Consider $m = \max\{\fidim(A), \fidim{B}\}$. Let $\mathcal{A}$ and $\mathcal{B}$ be the subgroups of $K_0(C)$ generated by the indecomposable modules of $\mod A$ and $\mod B$ respectively such that are not in $\mathcal{O}$. 
Note that $K_1(C) \subset \mathcal{A} \times \mathcal{B} \times  \mathcal{O}$, and
$$\bar{\Omega}_C ([M]) = \left\lbrace \begin{array}{ll} \bar{\Omega}_A([M]) + [N] & \text{ if } [M]\in \mathcal{A}, \\  \bar{\Omega}_B([M]) + [N]  & \text{ if } [M]\in \mathcal{B}\text{ and}  \\ {} [N] & \text{ if } [M]\in \mathcal{O}, \end{array} \right.$$
with $[N] \in \mathcal{O}$.

Consider $v=(v_1,v_2,v_3) \in K_1(C) \subset \mathcal{A} \times \mathcal{B} \times \mathcal{O}$ such that $\bar{\Omega}_C^n(v) = (0,0,0)$. Since $\fidim(A)<\infty$ and $\fidim(B)<\infty$, then $\bar{\Omega}_C^k(v) = (0,0,v')$ for $k\geq m$. On the other hand $\bar{\Omega}_C^{n-m}(0,0,v') = (0,0,0)$, then, by Proposition \ref{invariante} item 2, $n-m \leq \sup\{\eta_{\bar{\Omega}_C}(D):D\in \mathcal{O}\}$. We conclude that $\fidim(C) \leq \sup\{\eta_{\bar{\Omega}_C}(D):D\in \mathcal{O}\}+m+1$. Since $\mathcal{O}$ is syzygy finite, then $\fidim(C)< \infty$.\end{proof}

\end{teo}

As a particular case of the previous result we have.

\begin{coro}

Let $A = \frac{\Bbbk Q_A}{I_A}$ and $B = \frac{\Bbbk Q_B}{I_B}$ be finite dimensional algebras. Consider $C = \frac{\Bbbk \Gamma}{I_C}$ with the following conditions:

\begin{itemize}

\item ${Q_C}_0 = {Q_A}_0 \cup {Q_B}_0$.

\item ${Q_C}_1 = {Q_A}_1 \cup {Q_B}_1 \cup \{\alpha_j : \start(\alpha_j) \in  {Q_A}_0, \ \target(\alpha_j) \in {Q_B}_0 \}_{j \in J} \cup  \{\beta_k : \start(\beta_k) \in {Q_B}_0, \ \target(\beta_k) \in {Q_A}_0 \}_{k \in K}$.

\item $\langle I_A, I_B, \alpha \alpha_j, \beta \beta_k \text{ for }\alpha \in {Q_A}_1, \beta\in {Q_B}_1, \alpha_j\beta_k, \beta_k\alpha_j \text{ where } j\in J, k\in K\rangle \subset I_C$.

\item  $\Omega(S_v) = M_A \oplus N_B$ with $\pd_B(N_B) < \infty$ if $v\in {Q_A}_0$ and $\pd_A(M_A)<\infty$ if $v\in {Q_B}_0$.  

\item $\fidim(A) < \infty$ and $\fidim(B) < \infty$.

\end{itemize}

Then $\fidim(C) < \infty$.
\begin{proof}

Since $\Pi_A(\Omega_C(C_0))$ has finite projective dimension as $A$-module, then $[\add \Orb_{\Omega_A}(\Pi_A(\Omega_C(C_0)))]$ is a finitely generated subgroup of $K_0(A)$. Analogously we have that $[\add \Orb_{\Omega_B}( \Pi_B(\Omega_C(C_0)))]$ is a finitely generated subgroup of $K_0(B)$. Finally we deduce that $\mathcal{O} = [\add \Orb_{\Omega_A}(\Pi_A(\Omega_C(C_0))) \times \add \Orb_{\Omega_B}( \Pi_B(\Omega_C(C_0)))]$ is a finitely generated subgroup of $K_0(C)$.\end{proof}
\end{coro}

\begin{prop} \label{explicito} Let $A = \frac{\Bbbk Q_A}{I_A}$ and $B = \frac{\Bbbk Q_B}{I_B}$ be finite dimensional algebras. Consider $C = \frac{\Bbbk \Gamma}{I_C}$ with the following conditions:

\begin{itemize}

\item ${Q_C}_0 = {Q_A}_0 \cup {Q_B}_0$.

\item ${Q_C}_1 = {Q_A}_1 \cup {Q_B}_1 \cup \{\alpha_j : \start(\alpha_j) \in  {Q_A}_0, \ \target(\alpha_j) \in {Q_B}_0 \}_{j \in J} \cup  \{\beta_k : \start(\beta_k) \in {Q_B}_0, \ \target(\beta_k) \in {Q_A}_0 \}_{k \in K}$.

\item $I_C = \langle I_A, I_B, \alpha \alpha_j, \beta \beta_k \text{ for }\alpha \in {Q_A}_1, \beta\in {Q_B}_1, \alpha_j\beta_k, \beta_k\alpha_j \text{ where } j\in J, k\in K\rangle$.

\end{itemize}

Then 
\begin{enumerate}
\item $\max\{\phi\dim(A), \phi\dim(B)\} \leq \fidim(C) \leq \max\{\phi\dim(A), \phi\dim(B)\} + \vert {Q_A}_0\vert + \vert {Q_B}_0\vert +1$.

\item $\max\{\phi\dim(A^{op}), \phi\dim(B^{op})\} \leq \fidim(C^{op}) \leq \max\{\phi\dim(A^{op}), \phi\dim(B^{op})\} + \vert {Q_{A^{op}}}_0\vert + \vert {Q_{B^{op}}}_0\vert +1$.
\end{enumerate}

\begin{proof}
\begin{enumerate}

\item We denote by  $m = \max\{\phi\dim(A), \phi\dim(B)\}$,  $\partial A = \{v\in {Q_A}_0: \exists \ \alpha \in \Gamma_1 \text{ with } \start(\alpha) = v \text{ and } \target(\alpha) \in {Q_B}_0\}$ and $\partial B = \{v\in {Q_B}_0: \exists \ \beta \in \Gamma_1 \text{ with } \start(\beta) = v \text{ and } \target(\beta) \in {Q_A}_0\}$. For a vertex $v\in Q_C$ ($Q_A$, $Q_B$), ${P_v}_C$ (${P_v}_A$, ${P_v}_B)$ is the projective $C$-module ($A$-module, $B$-module) associate to the vertex $v$. 

$\mathcal{P}_{\partial A}$ is the full subcategory of $\mod A$ generated by the projective $A$-modules $P_v$ with $v \in \partial A$ and $ \mathcal{P}_{\partial B}$ is the full subcategory of $\mod B$ generated by the projective $B$-modules $P_v$ with $v \in \partial B$. 
For $M \in \mod A$, we have  $\Omega_C(M)= \Omega_A(M)\oplus P_B$. Analogously, for $M \in \mod B$ we have $\Omega_C(M)= \Omega_A(M)\oplus P_B$.

Note that we have $$K_1(C) \subset \langle [M]: M \in \mod A, [N]: N \in \mod B \rangle = [P_{\partial A}]\times [P_{\partial B}]\times K_0(A) \times K_0(B).$$

Let $v= (v_1, v_2, v_3, v_4) \in K_{1}(C)$ such that $\bar{\Omega}_C^{k-1}(v) \neq 0$ and $\bar{\Omega}_C^k(v) = 0$. Since $\bar{\Omega}_C^{\phi \dim (A)}(v) = (a,b,0,c)$ and $\bar{\Omega}_C^{\phi \dim (B)}(v) = (a',b',c',0)$, then $\bar{\Omega}_C^{m}(v)=(a'',b'',0,0)$. Since we have $\Omega_C(P_{\partial A}\oplus P_{\partial B}) \subset P_{\partial A} \oplus P_{\partial B}$, then, by Proposition \ref{invariante}, $\Omega_C^n(a'',b'',0,0) = 0$ where $n = \vert {Q_A}_0\vert + \vert {Q_B}_0\vert$. Therefore $k \leq \max\{\phi\dim(A), \phi\dim(B)\}+n$.

Suppose now that $\fidim(A) = m$, $\bar{\Omega}_A^m(v) = 0$ and $\bar{\Omega}_A^{m-1}(v) \not=0$ with $v\in K_0(A)$. Then we have $\bar{\Omega}^{m-1}_C(w_A,w_B,v,0)\not=0$ and $\bar{\Omega}^m_C(w_A,w_B,v,0)\not=(w'_A,w'_B,0,0)$ with $w_A \in [P_{\partial A}]$ and $w_B \in [P_{\partial B}]$.
We have two cases

\begin{itemize}

\item There exists $k \in \mathbb{N}$ such that $\bar{\Omega}^k_C(w'_A,w'_B,0,0) = (0,0,0,0)$.

\item  There exists $k \in \mathbb{N}$ such that $\bar{\Omega}^k_C(w'_A,w'_B,0,0) = (w''_A,w''_B,0,0)$ where $\forall s \in \mathbb{N}$ there is $({w_s}_A,{w_s}_B,0,0)$ such that $\bar{\Omega}_C^s({w_s}_A,{w_s}_B,0,0) = (w''_A,w''_B,0,0)$.

\end{itemize}

The first case implies $m \leq \fidim(C)$.
For the second case we can consider $(w_A-{w_{n+k}}_A,w_B-{w_{n+k}}_B,v,0) \in [P_{\partial A}]\times [P_{\partial B}]\times K_0(A) \times K_0(B)$, then $\bar{\Omega}^{n+k}(w_A-{w_{n+k}}_A,w_B-{w_{n+k}}_B,v,0)=(0,0,0,0)$, and $\bar{\Omega}^{n}(w_A-{w_{n+k}}_A,w_B-{w_{n+k}}_B,v,0)\not=(0,0,0,0)$ so we conclude that $m =\fidim(A) \leq \fidim(C)$. Analogously we conclude that $\fidim(B)\leq \fidim(C)$.\\

\item From now on we denote the arrows of $Q_{C^{op}}$ with the same name as $Q_{C}$. 

\underline{Claim 1:} If $M = \Omega (N)$ and $M = (M^v, T_{\gamma}) \in \ind C$, then $\topp M \in \mod A$ or $\topp M \in \mod B$.\\

Suppose $\topp M \not\in \mod A$ and $\topp M \not\in \mod B$. This implies that there exist vertices $v_0 \in {Q_C}_0$, $v_1 \in {Q_A}_0$, $v_2 \in {Q_B}_0$ and paths $\rho_1: v_1 \rightarrow v_0 $ and $\rho_2:v_2 \rightarrow v_0$ such that $T_{\rho_i}(x_i) = x_0 \neq 0$ for $i = 1,2$ and $x_i \in M^{v_i}$ for $i = 0,1,2$. Consider an indecomposable projective module $P = (P^v, \bar{T}_{\gamma})$ and a non trivial morphism $f: M \rightarrow P$ such that $f_{v_0}(x_0) = y_0$ with $y_0 \in P^{v_0}$. We conclude that there are non null vectors $y_1 \in P^{v_1}$ and $y_2 \in P^{v_2}$ such $T_{\rho_i}(y_i) = y_0$ for $i=1,2$, which is a contradiction since $\topp P$ is a simple module.\\

%

\underline{Claim 2:} Let $P = (P^v, T_{\gamma})\in \mathcal{P}_C$ be indecomposable with $\topp P \in \mod A  $. If $v_0 \in {Q_{A}}_0 $ is such that $P^{v_0} \not = 0$, then $T_{\beta_j}$ is injective for all $\beta_j:v_0 \rightarrow {Q_{B}}_0$.\\

It follows since for every nontrivial path $\rho:v_0 \rightarrow w$ in $Q_A$, $\rho \beta_j\not =0$ whenever $\target(\rho) = \start(\beta_j)$. 

Now consider $M$ non simple indescomposable module such that $M = \Omega(N)$ with $\topp(M) \in \mod A$. Then $M$ is a submodule of a projective module $P$ with $\topp P \subset \mod A$ and we conclude that claim 2 is also true for the module $M$.\\

We denote by $\mathcal{T} = \{ S_{v_0}\in \mathcal{S}(C) \text{ and there is } \beta_k \text{ or } \alpha_j \text{ with } \start(\beta_k) = v_0 \text{ or } \start(\alpha_j) = v_0 \}$.
Let $f: \Omega(\mod C^{op}) \rightarrow \mod A^{op} \times \mod B^{op} \times \langle C_0 \rangle$ be the function defined in the indecomposable modules as follows
$$f(M) = \left\lbrace \begin{array}{ll} (0,0,S_{v_0}) & \text{ if } M = S_{v_0}\in \mathcal{T} \\ (\prod_A(M),0,0) & \text{ if we are not in the first case and } \topp(M) \in \mod A \\ (0, \prod_B(M),0) & \text{ if we are not in the first case and } \topp(M) \in \mod B \end{array} \right.$$

\underline{Claim 3:} $f:  \Omega(\mod C^{op}) \rightarrow \mod A^{op} \times \mod B^{op} \times \langle C_0 \rangle$ is a monomorphism of groups and $f(\Omega_{C^{op}}(M)) = (\Omega_{A^{op}}(M_1), \Omega_{B^{op}}(M_2), 0)$ if $f(M) = (M_1, M_2, 0)$. 

If $M = M_1\oplus M_2 \oplus S$ where 

\begin{itemize}
\item $S$ is a direct sum of modules of $\mathcal{T}$.

\item $M_1$ has no direct summand that are simple modules of $\mathcal{T}$ and $\topp(M_1) \in \mod A^{op}$.

\item $M_2$ has no direct summand that are simple modules of $\mathcal{T}$ and $\topp(M_1) \in \mod B^{op}$.
\end{itemize}  

Then we have

\begin{itemize}
\item Let $N_1 = \Omega_{C^{op}}(M_1)$, then it has not direct summand that are simple modules of $\mathcal{T}$ and $\topp (N_1) \subset \mod A^{op}$ and the following diagram is commutative with exact rows.

$$\xymatrix{0 \ar[r] & \Omega_{A^{op}}(\prod_A(N_1)) \ar[r] &\prod_{A^{op}}(P) \ar[r] & \prod_{A^{op}}(M_1) \ar[r] & 0 \\ 0 \ar[r] & N_1 \ar[r]\ar[u]^{\prod_{A^{op}}} & P \ar[r] \ar[u]^{\prod_{A^{op}}} & M_1 \ar[r] \ar[u]^{\prod_A^{op}}&0 }$$ 

\item Let $N_2 = \Omega_{C^{op}}(M_2)$, then it has not direct summand that are simple modules of $\mathcal{T}$ and $\topp (N_2) \subset \mod {B^{op}}$. and the following diagram is commutative with exact rows.

$$\xymatrix{0 \ar[r] & \Omega_{B^{op}}(\prod_{B^{op}}(N_2)) \ar[r] &\prod_{B^{op}}(P) \ar[r] & \prod_{B^{op}}(M_2) \ar[r] & 0 \\ 0 \ar[r] & N_2 \ar[r]\ar[u]^{\prod_{B^{op}}} & P \ar[r] \ar[u]^{\prod_{B^{op}}} & M_2 \ar[r] \ar[u]^{\prod_{B^{op}}}&0 }$$ 

\end{itemize}
If $S=0$ then $f(\Omega_{C^{op}}(M)) = (\Omega_{A^{op}}(M_1), \Omega_{B^{op}}(M_2), 0)$.\\

On the other hand if $f(M_1\oplus M_2 \oplus S) = f(M'_1\oplus M'_2 \oplus S')$, then  $\prod_{A^{op}}(M_1) = \prod_{A^{op}}(M'_1)$ , $\prod_{B^{op}}(M_2) = \prod_{B^{op}}(M'_2)$ and $S = S'$. Since $ \prod_{A^{op}}(M_1) = \prod_{A^{op}}(M'_1)$ and $T_{\alpha_k}$ are injective for every representation in $\Omega(\mod {C^{op}})$, then $M_1 = M'_1$. The same follows from $\prod_{B^{op}}(M_2) = \prod_{B^{op}}(M'_2)$.

Finally we deduce, with similar computations to part 1, that $$\fidim(C^{op}) \leq \max\{\phi\dim(A^{op}), \phi\dim(B^{op})\} +  \vert {Q_{A^{op}}}_0\vert + \vert {Q_{B^{op}}}_0\vert + 1$$
\end{enumerate}

\end{proof}

\end{prop}

\begin{obs}\label{observacion} In the Hypothesis of Proposition \ref{explicito}, we have 

\begin{enumerate}
\item $\fidim(C) < \infty \Leftrightarrow \fidim(A) < \infty \text{ and } \fidim(B) < \infty.$

\item From the proof of Proposition \ref{explicito} we deduce that if $A$ and $B$ ($A^{op}$ and $B^{op}$) are syzygy finite then $C$ ($C^{op}$) is also syzygy finite.
\end{enumerate}
\end{obs}

We finish this section with two results on Morita context algebras for the finitistic dimenison.
The first one, the theorem below, has its hypotheses similar to Theorem \ref{teo}.

\begin{teo}
Let $A = \frac{\Bbbk Q_A}{I_A}$ and $B = \frac{\Bbbk Q_B}{I_B}$ be finite dimensional algebras. Consider $C = \frac{\Bbbk \Gamma}{I_C}$ with the following conditions:

\begin{itemize}

\item ${Q_C}_0 = {Q_A}_0 \cup {Q_B}_0$.

\item ${Q_C}_1 = {Q_A}_1 \cup {Q_B}_1 \cup \{\alpha_j : \start(\alpha_j) \in {Q_A}_0, \ \target(\alpha_j) \in {Q_B}_0 \}_{j \in J} \cup  \{\beta_k : \start(\beta_k) \in {Q_B}_0, \ \target(\beta_k) \in {Q_A}_0 \}_{k \in K}$.

\item $\langle I_A, I_B, \alpha \alpha_j, \beta \beta_k \text{ for }\alpha \in {Q_A}_1, \beta\in {Q_B}_1, \alpha_j\beta_k, \beta_k\alpha_j \text{ where } j\in J, k\in K\rangle \subset I_C$.

\item $\mathcal{O}= [\add \Orb_{\Omega_A}(\Pi_A(\Omega_C(B_0))) \times \add \Orb_{\Omega_B}( \Pi_B(\Omega_C(A_0)))] \subset K_0(C)$ is syzygy finite.


\item $\findim(A) < \infty$ and $\findim(B) < \infty$.

\end{itemize}

Then $\findim(C) < \infty$.

\begin{proof}

We begin the proof with a straightforward claim, it follows from Lemma \ref{lemita}.\\ 

\underline{Claim:} If $M \in \mod A$ ($\mod B$), then $\Omega_C(M) = \Omega_A(M) \oplus N$ ($ \Omega_B(M) \oplus N$) where $[N]\in \mathcal{O}$.\\

Consider $M \in \mod A$ ($\mod B$) such that $\pd_C(M)<\infty$. Since $\pd_A(M) \leq \pd_C(M)+ \findim\{N : [N] \in \mathcal{O}\}$ and $\mathcal{O}$ is syzygy finite, then $\pd_A (M) < \infty$ ($\pd_B(M)  < \infty$).   

Let $M\in \mod C$ such that $\pd_C(M) < \infty$. 
By Lemma \ref{lemita} $\Omega_C (M)= N_A \oplus N_B$ where $N_A\in \mod A$, $N_B\in \mod B$ and $\Omega^k_A (N_A) = \Omega^k_B(N_B) = 0$ if $k = \max\{\findim(A), \findim(B)\}$.

Finally, by the claim above $[\Omega^{k+1}(M)] \in \mathcal{O}$. Hence $\pd (\Omega^{k+1}(M)) < \infty$ and in particular $$\pd(M) \leq k+1 + \findim\{N : [N] \in \mathcal{O}\}.$$
\end{proof}
\end{teo}

The following proposition has its hypotheses similar to Theorem \ref{explicito}.

\begin{prop}
Let $A = \frac{\Bbbk Q_A}{I_A}$ and $B = \frac{\Bbbk Q_B}{I_B}$ be finite dimensional algebras. Consider $C = \frac{\Bbbk \Gamma}{I_C}$ with the following conditions:

\begin{itemize}

\item ${Q_C}_0 = {Q_A}_0 \cup {Q_B}_0$.

\item ${Q_C}_1 = {Q_A}_1 \cup {Q_B}_1 \cup \{\alpha_j : \start(\alpha_j) \in {Q_A}_0, \ \target(\alpha_j) \in {Q_B}_0 \}_{j \in J} \cup  \{\beta_k : \start(\beta_k) \in {Q_B}_0, \ \target(\beta_k) \in {Q_A}_0 \}_{k \in K}$.

\item $\langle I_A, I_B, \alpha \alpha_j, \beta \beta_k \text{ for }\alpha \in {Q_A}_1, \beta\in {Q_B}_1, \alpha_j\beta_k, \beta_k\alpha_j \text{ where } j\in J, k\in K\rangle = I_C$.

\item $\mathcal{O}= [\add \Orb_{\Omega_A}(\Pi_A(\Omega_C(B_0))) \times \add \Orb_{\Omega_B}( \Pi_B(\Omega_C(A_0)))] \subset K_0(C)$ is syzygy finite.


\item $\findim(A^{op}) < \infty$ and $\findim(B^{op}) < \infty$.

\end{itemize}

Then $\findim(C^{op}) < \infty$.

\begin{proof}

Consider $M \in \Omega(\mod C^{op})$, such that $\pd_{C^{op}} (M) < \infty$.
Note, by Claim 3 of Proposition \ref{explicito}, that for $M \in \Omega(\mod C^{op})$ we have the following:

\begin{itemize}
\item $M = M_1 \oplus M_2 \oplus S$ where $M_1 \in \mod A^{op}$, $M_2 \in \mod B^{op}$, and $S$ is a semisimple module, and 

\item  $\Omega_{C^op}^k(f(M_1)) = f(\Omega_{A^{op}}^k(M_2))$, $\Omega_{C^op}^k(f(M_2)) = f(\Omega_{B^{op}}^k(M_2))$.
\end{itemize}

Since $\pd_{C^{op}} (M) < \infty$, then $\pd_{C^{op}} (M_1) < \infty$ and $\pd_{C^{op}} (M_2) < \infty$, hence $\pd_{C^{op}} (M_1) < \pd_{A^{op}}(M_1) < \findim (A^{op})$ and $\pd_{C^{op}} (M_2) < \pd_{B^{op}}(M_2) < \findim (B^{op})$. Since $S \in \add C^{op}_0$ and it has finite nonisomorphic indecomposable modules, then $\pd_{C^{op}} (M) < \max\{\findim(A^{op}), \findim(B^{op}), \findim(C^{op}_0)\}$ and we conclude that $\findim(C^{op}) < \max\{\findim(A^{op}), \findim(B^{op}), \findim(C_0)\} +1$

\end{proof}

\end{prop}

\section{(Counter)examples}

In this section we present two examples. The first one (Example 5.4 of \cite{BM1}) shows that Hypotheses {\bf H1}, {\bf H2} and {\bf H4} are not enough to obtain the results of the previous section for the Igusa-Todorov $\phi$ function. The second is an example of an algebra where its $\phi$-dimension is not symmetric, i.e. $\fidim(A)$ and $\fidim(A^{op})$ do not agree.

\begin{ej} Let $C = \frac{\Bbbk Q}{I}$ be an algebra where $Q$ is 

$$
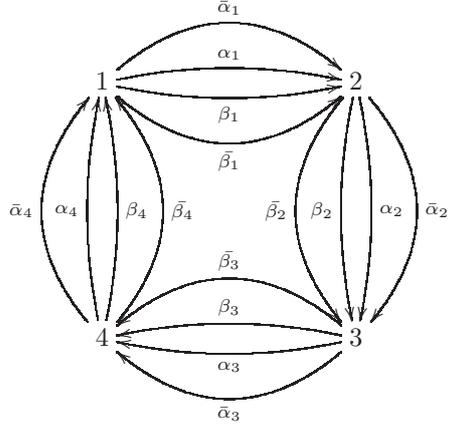
\xymatrix{ 1 \ar@/^8mm/[rrr]^{\bar{\alpha}_1} \ar@/^2mm/[rrr]^{\alpha_1} \ar@/_2mm/[rrr]_{\beta_1} \ar@/_8mm/[rrr]_{\bar{\beta_1}} & &  & 2 \ar@/^8mm/[ddd]^{\bar{\alpha}_2,} \ar@/^2mm/[ddd]^{\alpha_2} \ar@/_2mm/[ddd]_{\beta_2} \ar@/_8mm/[ddd]_{\bar{\beta_2}} \\ & &  &\\& & & & \\ 4 \ar@/^8mm/[uuu]^{\bar{\alpha}_4} \ar@/^2mm/[uuu]^{\alpha_4} \ar@/_2mm/[uuu]_{\beta_4} \ar@/_8mm/[uuu]_{\bar{\beta_4}} &  & &  3 \ar@/^8mm/[lll]^{\bar{\alpha}_3} \ar@/^2mm/[lll]^{\alpha_3} \ar@/_2mm/[lll]_{\beta_3} \ar@/_8mm/[lll]_{\bar{\beta_3}}}$$

and $I = \langle \alpha_{i}\alpha_{i+1}-\bar{\alpha}_{i}\bar{\alpha}_{i+1},\ \beta_{i}\beta_{i+1}-\bar{\beta}_{i}\bar{\beta}_{i+1},\ \alpha_{i}\bar{\alpha}_{i+1},\ \bar{\alpha}_{i}\alpha_{i+1},\ \beta_{i}\bar{\beta}_{i+1},\ \bar{\beta}_{i}\beta_{i+1},\text{ for } i \in \mathbb{Z}_4, \ J^3 \rangle$

Consider $A = \Bbbk Q_A$ and $B = \Bbbk Q_B$ where 
\begin{itemize}
\item $Q_A$ is the full quiver formed by the vertices $1$ and $2$.
\item $Q_B$ is the full quiver formed by the vertices $3$ and $4$.
\end{itemize}
It is easy to show that hypotheses {\bf H1}, {\bf H2} and {\bf H4} are verified. However it was showed in \cite{BM1} that $\fidim (C) = \infty$.

\end{ej}

\begin{ej}Consider the finite dimensional algebra $C_{p,q} =\frac{\Bbbk Q}{I_{p,q}}$ where $Q$ is the following quiver

$$
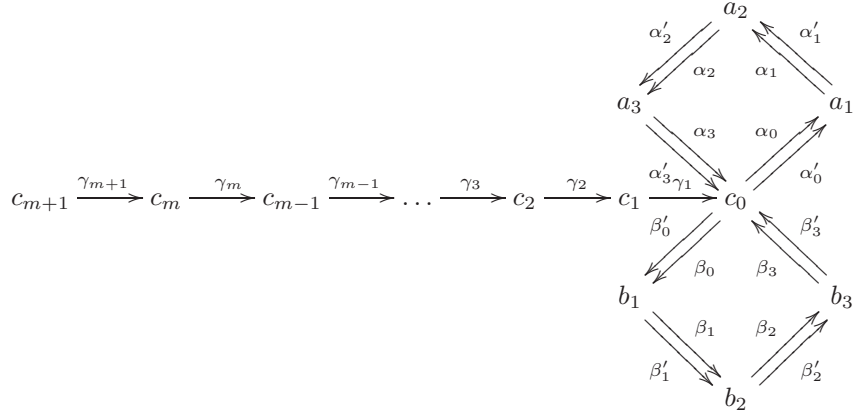
\xymatrix{ &&&&&& a_2 \ar@<.5ex>[dl]^{\alpha_2} \ar@<-.5ex>[dl]_{\alpha'_2} && \\&&&&& a_3 \ar@<.5ex>[dr]^{\alpha_3} \ar@<-.5ex>[dr]_{\alpha'_3} && a_1 \ar@<.5ex>[ul]^{\alpha_1} \ar@<-.5ex>[ul]_{\alpha'_1}& \\ c_{m+1} \ar[r]^{\gamma_{m+1}} &  c_m \ar[r]^{\gamma_m}&c_{m-1}\ar[r]^{\gamma_{m-1}}& \ldots \ar[r]^{\gamma_3}& c_2 \ar[r]^{\gamma_2}& c_1 \ar[r]^{\gamma_1}& c_0\ar@<.5ex>[ur]^{\alpha_0} \ar@<-.5ex>[ur]_{\alpha'_0} \ar@<.5ex>[dl]^{\beta_0} \ar@<-.5ex>[dl]_{\beta'_0} &&  \\&&&&& b_1 \ar@<.5ex>[dr]^{\beta_1} \ar@<-.5ex>[dr]_{\beta'_1} && b_3 \ar@<.5ex>[ul]^{\beta_3} \ar@<-.5ex>[ul]_{\beta'_3}& \\&&&&&& b_2 \ar@<.5ex>[ur]^{\beta_2} \ar@<-.5ex>[ur]_{\beta'_2} && }$$


and $I_{p,q}$ is generated by

\begin{itemize}
\item $\gamma_{i+1}\gamma_{i}\ \forall i=1,\ldots,m$,

\item $\alpha_{i}\alpha_{i+1}, \beta_{i}\beta_{i+1}, \alpha'_{i}\alpha'_{i+1}, \beta'_{i}\beta'_{i+1}\ \forall i \in \mathbb{Z}_4$,

\item $\alpha'_{i}\alpha_{i+1} - \alpha_{i}\alpha'_{i+1}, \beta'_{i}\beta_{i+1}-\beta_{i}\beta'_{i+1}\ \forall i \in \mathbb{Z}_4$,

\item $ \alpha_3\beta_0 , \alpha'_3\beta'_0, \beta_3\alpha_0, \beta'_3\alpha'_0$ and
\item $\alpha'_3\beta_0 -\alpha_3\beta'_0, \beta'_3\alpha_0-p\beta_3\alpha'_0, \gamma_1\alpha_0-q\gamma_1\alpha'_0, \gamma_1\beta'_0-\gamma_1\beta_0$ with $\mathbb{Q}[p,q] \cong \mathbb{Q}[x,y]$. 
\end{itemize}

Since $\rad^3 (C_{p,q}) = 0$, then $\ll(M) \leq 3$, and $\ll(\Omega(M)) \leq 2$ for all $ M \in \mod C_{p,q}$.\\

Let $\bar{Q} \subset Q$ be the full subquiver with vertices $\{ a_1, a_2, a_3, b_1, b_2, b_3, c_0\}$. Consider $B_{p,q} = \frac{\Bbbk \bar{Q}}{I_{p,q} \cap \Bbbk \bar{Q}}$. 
If $M$ is in $\mod B$, then the $\topp$ of the indecomposable modules of $\Omega (M)$ is in

\begin{itemize}

\item $c_0$ (level 0),

\item $a_1$ or $b_1$ (level 1),

\item $a_2$ or $b_2$ (level 2) and

\item $a_3$ or $b_3$ (level 3).

\end{itemize}

Consider the following modules 

$$M_{i,\lambda, n} \text{ such that } \supp(M_{1,\lambda, n}) = \xymatrix{a_i \ar@<.5ex>[r]^{\alpha_i}  \ar@<-.5ex>[r]_{\alpha'_i} & a_{i+1}} \text{ and } M_{i,\lambda, n} = \xymatrix{ \cdots \Bbbk^n \ar@<.5ex>[r]^{J_{\lambda}} \ar@<-.5ex>[r]_{1_{\Bbbk^n}} & \Bbbk^n \cdots}$$

$$M'_{i,\lambda, n} \text{ such that } \supp(M'_{1,\lambda, n}) = \xymatrix{a_i \ar@<.5ex>[r]^{\alpha_i}  \ar@<-.5ex>[r]_{\alpha'_i} & a_{i+1}} \text{ and } M'_{i,\lambda, n} = \xymatrix{ \cdots \Bbbk^n \ar@<.5ex>[r]^{{1_{\Bbbk^n}}} \ar@<-.5ex>[r]_{J_{\lambda}} & \Bbbk^n \cdots}$$

$$M_{i, n} \text{ such that } \supp(M_{1, n}) = \xymatrix{a_i \ar@<.5ex>[r]^{\alpha_i}  \ar@<-.5ex>[r]_{\alpha'_i} & a_{i+1}} \text{ and } M_{i, n} = \xymatrix{ \cdots \Bbbk^n \ar@<.5ex>[r]^{(1_{{\Bbbk}^n},0)} \ar@<-.5ex>[r]_{(0, 1_{{\Bbbk}^n})} & \Bbbk^{n+1} \cdots}$$

$$\bar{M}_{i, n} \text{ such that } \supp(\bar{M}_{1, n}) = \xymatrix{a_i \ar@<.5ex>[r]^{\alpha_i}  \ar@<-.5ex>[r]_{\alpha'_i} & a_{i+1}} \text{ and } \bar{M}_{i, n} = \xymatrix{ \cdots \Bbbk^{n+1} \ar@<.5ex>[r]^{(1_{{\Bbbk}^n},0)^t} \ar@<-.5ex>[r]_{(0, 1_{{\Bbbk}^n})^t} & \Bbbk^{n}\cdots}$$

$$N_{i,\lambda, n} \text{ such that } \supp(N_{1,\lambda, n}) = \xymatrix{b_i \ar@<.5ex>[r]^{\beta}  \ar@<-.5ex>[r]_{\beta'_i} & b_{i+1}} \text{ and } N_{i,\lambda, n} = \xymatrix{ \cdots \Bbbk^n \ar@<.5ex>[r]^{J_{\lambda}} \ar@<-.5ex>[r]_{1_{\Bbbk^n}} & \Bbbk^n \cdots}$$

$$N'_{i,\lambda, n} \text{ such that } \supp(N'_{1,\lambda, n}) = \xymatrix{b_i \ar@<.5ex>[r]^{\beta}  \ar@<-.5ex>[r]_{\beta'_i} & b_{i+1}} \text{ and } N'_{i,\lambda, n} = \xymatrix{ \cdots \Bbbk^n \ar@<.5ex>[r]^{1_{\Bbbk^n}} \ar@<-.5ex>[r]_{J_{\lambda}} & \Bbbk^n \cdots}$$

$$N_{i, n} \text{ such that } \supp(N_{1, n}) = \xymatrix{b_i \ar@<.5ex>[r]^{\beta_i}  \ar@<-.5ex>[r]_{\beta'_i} & b_{i+1}} \text{ and } N_{i, n} = \xymatrix{ \cdots \Bbbk^n \ar@<.5ex>[r]^{(1_{{\Bbbk}^n},0)} \ar@<-.5ex>[r]_{(0, 1_{{\Bbbk}^n})}  & \Bbbk^{n+1}\cdots}$$

$$\bar{N}_{i, n} \text{ such that } \supp(\bar{N}_{1, n}) = \xymatrix{b_i \ar@<.5ex>[r]^{\beta_i}  \ar@<-.5ex>[r]_{\beta'_i} & b_{i+1}} \text{ and } \bar{N}_{i, n} = \xymatrix{ \cdots \Bbbk^{n+1} \ar@<.5ex>[r]^{(1_{{\Bbbk}^n},0)^t} \ar@<-.5ex>[r]_{(0, 1_{{\Bbbk}^n})^t} & \Bbbk^{n} \cdots} $$
for $i = 1,2,3$, $\lambda \in \Bbbk$, $n \in \mathbb{N}$ and $a_4=b_4=c_0$ 

$$M_{0,\lambda, \mu, n} \text{ such that } \supp(M_{0,\lambda, \mu, n}) = \xymatrix{ b_1 &  c_0 \ar@<.5ex>[l]^{\beta'_0}  \ar@<-.5ex>[l]_{\beta_0} \ar@<.5ex>[r]^{\alpha_0}  \ar@<-.5ex>[r]_{\alpha'_0} & a_{1}} \text{ and } M_{0,\lambda, \mu, n} = \xymatrix{ \cdots \Bbbk^n & \Bbbk^n  \ar@<.5ex>[l]^{1_{\Bbbk^n}} \ar@<-.5ex>[l]_{J_{\lambda}} \ar@<.5ex>[r]^{J_{\mu}} \ar@<-.5ex>[r]_{1_{\Bbbk^n}} & \Bbbk^n \cdots}$$

$$M'_{0,\lambda, \mu, n} \text{ such that } \supp(M'_{0,\lambda, \mu, n}) = \xymatrix{ b_1 &  c_0 \ar@<.5ex>[l]^{\beta'_0}  \ar@<-.5ex>[l]_{\beta_0} \ar@<.5ex>[r]^{\alpha_0}  \ar@<-.5ex>[r]_{\alpha'_0} & a_{1}} \text{ and } M'_{0,\lambda, \mu, n} = \xymatrix{ \cdots \Bbbk^n & \Bbbk^n  \ar@<.5ex>[l]^{J_{\lambda}} \ar@<-.5ex>[l]_{1_{\Bbbk^n}} \ar@<.5ex>[r]^{1_{\Bbbk^n}} \ar@<-.5ex>[r]_{J_{\mu}} & \Bbbk^n \cdots}$$

$$M_{0, n} \text{ such that } \supp(M_{0, n}) = \xymatrix{ b_1 &  c_0 \ar@<.5ex>[l]^{\beta'_0}  \ar@<-.5ex>[l]_{\beta_0} \ar@<.5ex>[r]^{\alpha_0}  \ar@<-.5ex>[r]_{\alpha'_0} & a_{1}} \text{ and } M_{0, n} = \xymatrix{ \cdots \Bbbk^{n+1} & \Bbbk^{n} \ar@<.5ex>[l]^{(0,1_{{\Bbbk}^n})} \ar@<-.5ex>[l]_{(1_{{\Bbbk}^n},0)}\ar@<.5ex>[r]^{(1_{{\Bbbk}^n},0)} \ar@<-.5ex>[r]_{(0, 1_{{\Bbbk}^n})} & \Bbbk^{n+1} \cdots}$$

$$\bar{M}_{0, n} \text{ such that } \supp(M_{0, n}) = \xymatrix{ b_1 &  c_0 \ar@<.5ex>[l]^{\beta'_0}  \ar@<-.5ex>[l]_{\beta_0} \ar@<.5ex>[r]^{\alpha_0}  \ar@<-.5ex>[r]_{\alpha'_0} & a_{1}} \text{ and } \bar{M}_{0, n} = \xymatrix{ \cdots \Bbbk^{n} & \Bbbk^{n+1} \ar@<.5ex>[l]^{(0,1_{{\Bbbk}^n})^t} \ar@<-.5ex>[l]_{(1_{{\Bbbk}^n},0)^t}\ar@<.5ex>[r]^{(1_{{\Bbbk}^n},0)^t} \ar@<-.5ex>[r]_{(0, 1_{{\Bbbk}^n})^t} & \Bbbk^{n}\cdots}$$
for $\lambda \in \Bbbk$ and $n \in \mathbb{N}$.\\

Note that 

\begin{itemize}

\item $M_{1,0} = \bar{M}_{2,0} = S_{a_2}$, $M_{2,0} = \bar{M}_{3,0} = S_{a_3}$, $M_{3,0} = \bar{M}_{0,0} = S_{c_0}$ 

\item $N_{1,0} = \bar{N}_{2,0} = S_{b_2}$, $N_{2,0} = \bar{N}_{3,0} = S_{b_3}$, $N_{3,0} = \bar{M}_{0,0} = S_{c_0}$, $M_{0,0} = \bar{M}_{1,0}\oplus \bar{N}_{1,0} = S_{a_1}\oplus S_{b_1}$ 

\item $M'_{i,\lambda,n} = M_{i,\frac{1}{\lambda},n}, N'_{i,\lambda,n} = N_{i,\frac{1}{\lambda},n}, M'_{0,\lambda, \mu, n}= M_{0, \frac{1}{\lambda}, \frac{1}{\mu}, n}$ if $\mu \not = 0, \lambda \not = 0$.

\item $\Omega(M_{1, \lambda, n}) = M_{2, \lambda, n}$, $\Omega(M_{2, \lambda, n}) = M_{3, \lambda, n}$, $\Omega(M_{3, \lambda, n}) = M_{0,\lambda, \lambda, n}$ and $\Omega(M_{0,\lambda, \mu, n}) = M_{1, \lambda, n} \oplus N_{1, \mu, n}$.

\item $\Omega(M'_{1, 0, n}) = M'_{2, 0, n}$, $\Omega(M'_{2, 0, n}) = M'_{3, 0, n}$, $\Omega(M'_{3, 0, n}) = M'_{0, 0, 0, n}$ and $\Omega(M'_{0, 0, 0, n}) = M'_{1, 0, n} \oplus N'_{1, 0, n}$.

\item $\Omega(M_{1,n}) = M_{2,n-1}$, $\Omega(M_{2,n}) = M_{3,n-1}$, $\Omega(M_{3, n}) = M_{0,n-1}$ and $\Omega(M_0,n) = M_{1,n-1} \oplus N_{1,n-1}$ for $n \geq 1$.

\item $\Omega(\bar{M}_{1,n}) = \bar{M}_{2,n+1}$, $\Omega(\bar{M}_{2,n}) = \bar{M}_{3,n+1}$, $\Omega(\bar{M}_{3, n}) = \bar{M}_{0,n+1}$ and $\Omega(\bar{M}_{0,n}) = \bar{M}_{1,n+1} \oplus \bar{N}_{1,n+1}$.

\item $\Omega(N_{1, \lambda, n}) = N_{2, \lambda, n}$, $\Omega(N_{2, \lambda, n}) = N_{3, \lambda, n}$ and $\Omega(N_{3, \lambda, n}) = M_{0, p\lambda, \lambda , n}$. 

\item $\Omega(N'_{1, 0, n}) = N'_{2, 0, n}$, $\Omega(N'_{2, 0, n}) = N'_{3, 0, n}$, and $\Omega(N'_{3, 0, n}) = M'_{0, 0, 0, n}$

\item $\Omega(N_{1,n}) = N_{2,n-1}$, $\Omega(N_{2,n}) = N_{3,n-1}$, $\Omega(N_{3, n}) = N_{0, n-1}$ and $\Omega(M_{0,n}) = M_{1,n-1} \oplus N_{1, n-1}$ for $n \geq 1$.

\item $\Omega(\bar{N}_{1,n}) = \bar{N}_{2,n+1}$, $\Omega(\bar{N}_{2,n}) = \bar{N}_{3,n+1}$, $\Omega(\bar{N}_{3, n}) = \bar{N}_{0,n+1}$ and $\Omega(\bar{N}_{0,n}) = \bar{M}_{1,n+1} \oplus \bar{N}_{1,n+1}$.

\end{itemize}

Cosider $\mathcal{V}$ the full subcategory generated by de direct sum of the following family of modules
$$\{M_{i,\lambda, n}, M'_{i, 0, n}, M_{i, n}, \bar{M}_{i, n}, N_{j,\lambda, n}, N'_{i, 0, n}, N_{j, n}, \bar{N}_{j, n}, M_{0, \lambda, \mu, n}, M'_{0, 0, 0, n} \text{ for } i = 0,1,2,3, j=1,2,3 \text{ and } \lambda \in \Bbbk \}$$

\underline{Claim 1:} $\fidim (\mathcal{V}) = 4$.\\

It follows from the fact that $\bar{\Omega}\vert_{\Omega^4(\mathcal{V})}$ is an monomorphism.\\

\underline{Claim 2:} $\Omega^2(\mod B_{p,q}) \subset V \subset \Omega (\mod C_{p,q})$ and $K_2(C_{p,q}) \subset \langle [\mathcal{V}], [S_{c_{m-1}}], [S_{c_{m-2}}], \ldots, [S_{c_1}]\rangle$\\

Note that $\Omega(S_{c_1}) = M_{0,q,1,1}$, then 
\begin{itemize}
\item $[\Omega(S_{c_1})] \not \in \bar{\Omega}([\mathcal{V}])$,
\item $[\Omega^2(S_{c_1})]  \not \in \langle \bar{\Omega}^2([\mathcal{V}]), [\Omega(S_{c_1})]\rangle$,
\item $[\Omega^3(S_{c_1})]  \not \in \langle \bar{\Omega}^3([\mathcal{V}]), [\Omega(S_{c_1})], [\Omega^2(S_{c_1})] \rangle$ and
\item $[\Omega^4(S_{c_1})]  \not \in \langle \bar{\Omega}^4([\mathcal{V}]), [\Omega(S_{c_1})], [\Omega^2(S_{c_1})],[\Omega^3(S_{c_1})] \rangle$.
\end{itemize}

\underline{Claim 3:} $\phi \dim (C_{p,q}) = 5$.\\

It follows from the fact that $\bar{\Omega} \vert_{\langle [\Omega(S_{c_1})], [\Omega^2(S_{c_1})], [\Omega^3(S_{c_1})], [\Omega^4(S_{c_1})], \bar{\Omega}^4 ([\mathcal{V}]) \rangle}$ is a monomorphism.\\

On the other hand $\fidim(C_{p,q}^{op}) \geq m-1$ because $\id_{C_{p,q}} (S_{c_1}) = m-1$.

\end{ej}

Note that $C_{p,q}$ also does not verify Proposition \ref{explicito}, just take $A = \frac{\Bbbk Q'}{J^2}$ where   
$Q'$ is the subquiver of $Q$ generated by $\{c_0, c_1, \ldots, c_{m+1}\}$ and $B_{p,q}$ as before.

\end{document}